\documentclass{article}
\usepackage{amssymb}
\title{\textsf{Multiplicative measures on free groups}}
\author{\textsf{Alexandre V. Borovik}\thanks{Supported
by the Royal Society Leverhulme Trust Senior Research
Fellowship.} \and \textsf{Alexei G. Myasnikov} \and
\textsf{Vladimir N. Remeslennikov\thanks{Supported by EPSRC grant
GR/R29451.}}}
\date{\textsf{ 11 April 2002}}
\newtheorem{theorem}{Theorem}[section]

\newtheorem{corollary}[theorem]{Corollary}
\newtheorem{lemma}[theorem]{Lemma}

\newcommand{\bnm}[2]{\left(\begin{array}{c}#1 \cr #2\end{array}\right)}
\begin{document}

\maketitle

\pagestyle{myheadings} \markright{{\scriptsize A.~V.~ Borovik,
A.~G.~Myasnikov and V.~N.~Remeslennikov $\bullet$ Multiplicative
measures  $\bullet$ 11.04.02 }}

\section{How one can measure subsets in the free group?}

\subsection{Motivation}

The present paper is motivated  by needs of
practical computations in finitely presented groups.
In particular, we wish to develop tools which can be used in the
analysis of the ``practical" complexity of algorithmic problems for discrete
infinite groups, as well as in the analysis of the
behaviour of heuristic (e.g. genetic)  algorithms for infinite groups
\cite{My3,My2}.

In most computer-based computations in finitely presented groups
$G = F/R$ the elements are represented as freely reduced words in
the free group $F$, with procedures for comparing their images in
the factor group  $G = F/R$. Therefore the ambient algebraic
structure in all our considerations is the free group  $F=F(X)$
on a finite set $X = \{x_1,\ldots, x_m\} $. We identify  $F$ with
the set of all freely reduced words in the alphabet $X \cup
X^{-1}$, with the multiplication given by concatenation of words
with the subsequent free reduction.

The most natural and convenient way to generate pseudorandom
elements in $G$ is to produce pseudorandom freely reduced words in
$F$. The most abstract mathematical model of a random word
generator in $F$ is just a probabilistic distribution on $F$.
We find ourselves in the setting of the paper \cite{measures},
which initiated a general discussion of
probabilistic measures on the free group.
Analysis of complexity of algorithms on groups  necessarily
involves the study of their behaviour with respect to the size of
the input, usually, the length of input words. Different
probabilistic distributions on $F$ represent pseudorandom
generators with varying mean length of words. This mean length is
one of the most important parameters of a pseudorandom generator.
Since we wish to vary the mean length of inputs, a single fixed
distribution on $F$ does not suffice, and  we need a parametric
family of probabilistic distributions ${\cal P} = \{P_l\}$  of
varying mean length $l$ of elements. This leads to the crucial
point of our approach: a  measure of a given subset $R \subset F$
is not a particular number $P_l(R)$ (which is usually
meaningless), but rather a function ${\cal P}(R): l \rightarrow
P_l(R)$ which naturally encodes all statistical properties of $R$
with respect to the family of distributions ${\cal P}$. It turns
out that such well-known asymptotic characteristics of $R$ as
asymptotic density, co-growth rate, etc., are just the standard
analytic characteristics of the function ${\cal P}$. This opens
the way to apply classical analytical methods for description of
statistical behaviour of algorithms in groups.
In Section~\ref{se:1-5} we introduce a hierarchy of subsets $R$ in $F$ with
respect to their size, which is based on linear approximations of
the function $\mu(R)$. This hierarchy is quite  sensitive, for
example, it allows one to differentiate between sets with the same
asymptotic density.

Our requirements to probabilistic distributions
are motivated by a very practical,
engineering approach to computations in groups.
First of all, the probabilistic distribution should not be unnatural
in the context of computational group theory. It should provide an easy way to
make crude estimates of probabilities of various subsets important in
standard problems of group theory: subgroups (first of all, normal or finitely
generated subgroups),
cosets with respect to subgroups,
conjugacy classes, sets of words of special nature
(say, squares or commutators).
It should also provide for an easy analysis
 of asymptotic behaviour of probabilities
 when the mean word length tends to infinity.

Many sets we wish to measure have happened to be  context
free languages \cite{rayward-smith}. An important subclass is made
of {\em regular} subsets (that is, subsets produced by
deterministic finite automata). This very natural class of subsets
includes finitely generated subgroups and their cosets, and {\em
finitely generated cones} (sets of all words which start with an
initial segment belonging to a given finite set of words). The
class of regular sets in $F$ is closed under Boolean operations,
and under translation and conjugation by elements of $F$.

\subsection{Generation of random words in $F$}
\label{se:1-2}

 Let $F = F(X)$ be a free group with basis $X = \{x_1, \ldots,
x_m\}$. We use, as our random word generator, the following
no-return random walk $W_s$ ($s \in (0,1]$) on the Cayley
graph $C(F,X)$ of $F$ with respect  to the generating set $X$.
 We start at the identity element
$1$ and either do nothing with probability
$s$ (and return value $1$ as the output of our random word generator),
or move to one of the $2m$ adjacent vertices with equal
probabilities $(1-s)/2m$. If we are at a vertex $v \ne 1$, we
either stop at $v$ with probability $s$ (and return $v$), or move, with
probability $\frac{1-s}{2m-1}$, to one of the $2m-1$ adjacent
vertices lying away from $1$, thus producing a new freely reduced
word $vx_i^{\pm 1}$. In other words, we make random freely
reduced words $w$ of random lengths $|w|$ distributed according
to the geometric law
$$
P(|w| = k) = s(1-s)^k,
$$
in such way that words of the same length $k$ are produced with
equal probabilities (in terminology of \cite{measures}, we say that our measure is {\em moderated\/} by the geometric distribution of $\mathbb{N} \cup \{\,0\,\}$).
Observe that the set of all  words of
length $k$ in $F$ forms the sphere $S_k$ of radius $k$ in
$C(F,X)$ of cardinality  $|S_k| = 2m(2m-1)^{k-1}$. It is easy to
see that the resulting probabilistic atomic
measure\footnote{Recall that a measure $\mu$ on a countable set
$X$ is {\em atomic} if every subset $Y \subseteq X$ is
measurable. This is equivalent to saying that every singleton
subset $\{x\}$ is measurable. Obviously, $\mu(Y)= \sum_{x\in Y}
\mu(x)$.}  $\mu_s$ on $F$ is given by the formula
\begin{equation}
\mu_s(w) = \frac{s(1-s)^{|w|}}{2m\cdot (2m-1)^{|w|-1}} \quad
\hbox{ for } w \ne 1
\end{equation}
and
\begin{equation}
\mu_s(1) =s.
\end{equation}

 Thus, $\mu_s(w)$ is the probability that the random walk $W_s$
stops at $w$.    The mean length $L_s$ of words in $F$
distributed according to $\mu_s$ is equal to $$L_s = \sum_{w \in
F} |w|\mu_s(w) = s \sum_{k= 1}^{\infty}k(1-s)^{k-1} = \frac{1}{s}
-1. $$ Hence we have a family of probabilistic distributions $\mu
= \{\mu_s\}$ with the stopping probability $s \in (0,1)$ as a
parameter, which is related to the  average length $L_s$ as
$$s = \frac{1}{L_s+1}.$$
By  $\mu(R)$ we denote the function
\begin{eqnarray*}
\mu(R): (0,1) & \rightarrow & {\mathbb{R}}\\
s & \mapsto & \mu_s(R);
\end{eqnarray*}
we call it  {\it measure} of $R$ with respect to the
family of distributions $\mu$.

Denote by $n_k = n_k(R) = |R \cap S_k|$ the number of elements
 of length $k$ in $R$, and by $f_k = f_k(R)$ the relative
frequencies
$$f_k = \frac{|R \cap S_k|}{|S_k|}$$ of words of length $k$ in
$R$. Notice that $f_0 = 1 $ or $0$ depending on whether $R$ contains $1$ or not.
Recalculating $\mu(R)$ in terms of $s$, we immediately come
to the formula $$ \mu(R) = s\sum_{k=0}^\infty f_k(1-s)^k, $$ and
the series on the right hand side is convergent for all $s \in
(0,1)$. Thus, for every subset $R \subseteq F$, $\mu(R)$ is an analytic function of $s$.
 When studying the behaviour of $\mu(R)$,
we mostly restrict it to real arguments, but occasionally need
to work with extensions of $\mu(R)$ to larger regions of
the complex plane.
We use only most basic facts of the theory of analytic functions
which can be found in any book on complex analysis.

Notice that the asymptotic behaviour of the set $R$ when $L_s
\rightarrow \infty$ corresponds to the behaviour of the function
$\mu(R)$ when $s\rightarrow 0^+$. This will be  discussed  in more
detail in Section \ref{se:1-4}.  Here we just mention how one can
obtain a first coarse  approximation of the asymptotic behaviour
of the function $\mu(R)$.   Let $W_0$ be the  no-return
non-stop simple random walk  on $C(F,X)$ (like $W_s$ with $s =
0$), where the walker moves from a given vertex  to any adjacent
vertex away from the initial vertex 1 with equal probabilities
$1/2m$. In this event, the probability $\lambda(w)$ that the
walker hits an element $w \in F$ in $|w|$ steps
(which is the same as the probability
that the walker ever hits $w$)  is equal to
$$
\lambda(w) = \frac{1}{2m(2m-1)^{|w|-1}}, \ \hbox{ if } \ w \neq 1, \ \
\hbox{ and }  \ \lambda(1) = 1.
$$
 This gives rise to an atomic measure
$$ \lambda(R) = \sum_{w \in R}\lambda(w) = \sum_{k=0}^{\infty}f_k(R)$$
where $\lambda(R)$  is just the sum of the relative frequencies
of $R$.
 This measure is not probabilistic,
since some sets have no finite measure (obviously, $\lambda(F) =
\infty$), moreover, the  measure $\lambda$ is finitely additive,
but not $\sigma$-additive. We shall call $\lambda$  the {\em
frequency } measure on $F$. If $R$ is $\lambda$-measurable (i.e.,
$\lambda(R) <  \infty$) then $f_k(R) \rightarrow 0$ when $k
\rightarrow \infty$, so intuitively, the set $R$ is "small" in
$F$.

 A number of papers (see, for example, \cite{Arzh}, \cite{Fr}, \cite{Olsh},
\cite{woess}), used the {\em asymptotic density} (or more,
precisely, the {\em spherical asymptotic density})
 $$\rho(R) = \lim\sup f_k(R)$$
  as a numeric characteristic of the set $R$ reflecting its asymptotic behavior.
  Unfortunately, the asymptotic density is not even finitely additive, and it is not
   sensitive enough: many interesting   sets have asymptotic density either $1$ or $0$.

 More subtle analysis of asymptotic
behaviour of $R$ in some cases  provides  the
 {\em relative growth  rate}
$$
\gamma(R) = \lim\sup \sqrt[k]{f_k(R)}.
$$
Notice the obvious inequality $\gamma(R) \leqslant 1$.
If
$\gamma(R) < 1$ (we will have to say more about this case in
Section \ref{sec:cogrowth}), then, by an elementary result from Calculus,
the series $\sum f_k$
converges. This shows that if  $\gamma(R) < 1$ then $R$ is
$\lambda$-measurable.

Our distribution $\mu_s$ has the uncomfortably big standard
deviation $\sigma = \frac{\sqrt{1-s}}{s}$. This reflects the fact
that it is strongly skewed towards `short' elements. However,
since real life computations take place in the vicinity of $1$,
we believe that our model is useful as a first step in developing
statistical approach to computational group theory.

\subsection{The multiplicativity of the measure and generating functions}
\label{se:1-3}

It is convenient to renormalise our measures $\mu_s \in \mu$  and
work with the parametric  family $\mu^* = \{\mu_s^* \}$ of  {\em
adjusted measures}
\begin{equation}
\mu_s^*(w) = \left(\frac{2m}{2m-1} \cdot \frac{1}{s} \right)\cdot
\mu_s(w).
\end{equation}
This new measure $\mu_s^*$  is {\it multiplicative} in the sense
that
\begin{equation}
\mu_s^*(u\circ v) = \mu_s^*(u)\mu_s^*(v),
\end{equation}
where $u\circ v$ denotes the product of non-empty words $u$ and
$v$ such that $|uv| = |u| +|v|$ i.e.\ there is no cancellation
between $u$ and $v$. The measure $\mu$ itself is {\em almost multiplicative} in the sense
that
\begin{equation}
\mu_s(u\circ v) = c\mu_s(u)\mu_s(v) \quad \hbox{for} \quad c = \frac{2m}{2m-1} \cdot \frac{1}{s}
\end{equation}
for all non-empty words $u$ and
$v$ such that $|uv| = |u| +|v|$.

If we denote
\begin{equation}
t = \mu_s^*(x_i^{\pm 1})=  \frac{1-s}{2m-1} \label{eq:adjusted}
\end{equation}
 then
\begin{equation}
\mu_s^*(w) = t^{|w|}
\end{equation}
 for every non-empty word $w$.

Similarly, we can adjust the frequency  measure $\lambda$
making it into a multiplicative atomic measure
\begin{equation}
\lambda^*(w) = \frac{1}{(2m-1)^{|w|}}. \label{eq:lambda*}
\end{equation}

Let now $R$ be a subset in $F$ and  $n_k = n_k(R) = |R \cap
S_k|$ be the number of elements of length $k$ in $R$.  The
sequence $\{n_k(R)\}_{k=0}^{\infty}$ is called the {\it spherical
growth sequence} of $R$. We assume, for the sake of minor
technical convenience, that $R$ does not contain the identity
element $1$, so that  $n_0 = 0$. It is easy to see now that $$
\mu^*(R) = \sum_{k=0}^\infty n_kt^k .$$  One can view
$\mu^*(R)$ as the generating function of the spherical growth
sequence of the set $R$  in variable $t$ which is convergent for
each  $t \in [0,1)$. This simple observation will allow us (see
Sections 3 and 4) to apply a well established machinery of
generating functions of context-free languages to estimate
probabilities  of sets.

\subsection{Cesaro density}
\label{se:1-4}

 Let $\mu = \{\mu_s\}$ be the parametric family of distributions defined
above. For a subset $R$ of  $F$ we define the {\it limit measure}
$\mu_0(R):$
 $$ \mu_0(R) = \lim_{s \rightarrow 0^+}\mu(R) =  \lim_{s \rightarrow 0^+} s \cdot \sum_{k=0}^\infty
f_k(1-s)^k.$$  The function $\mu_0$   is additive, but not
$\sigma$-additive, since $\mu_0(w) = 0$ for a single element $w$.
It is easy to construct a set $R$ such that $\lim_{s\rightarrow
0^+}\mu(R)$ does not exist. However, in the applications that we
have in mind we have not yet encountered such a situation.
Strictly speaking, $\mu_0$ is not a measure because the set  of
all $\mu_0$-measurable sets is not closed under intersections
(though it is closed under complements). Because $\mu_s(R)$ gives
an approximation of $\mu_0(R)$ when $s \rightarrow 0^+$, or
equivalently, when $L_s \rightarrow \infty$, we shall call $R$
{\it measurable at infinity} if $\mu_0(R)$ exists, otherwise $R$
is called {\it singular}.

If $\mu(R)$ can be expanded  as a convergent power series in $s$
at $s=0$ (and hence in some neighborhood of $s = 0$):
 $$ \mu(R) = m_0 + m_1s + m_2s^2 + \cdots, $$ then
  $$\mu_0(R) = \lim_{s
\rightarrow 0^+} \mu(R) = m_0.$$ A corollary from a theorem by Hardy
and Littlewood \cite[Theorem~94]{hardy} (see
Corollary~\ref{cor:Cesaro1} in Section \ref{se:6}) asserts that
$\mu_0$ can be computed as the {\em Cesaro limit}
 \begin{equation}
 \label{eq:cesaro-1}
  \mu_0(R)  = \lim_{n\rightarrow\infty} \frac{1}{n}\left(f_1+\cdots+f_n\right).
  \end{equation}
   So it will be also natural to call $\mu_0$ the {\em Cesaro density}, or
    {\it asymptotic average density}.

 Notice,
that the Cesaro density $\mu_0$ is  more sensitive then the
standard asymptotic density $\rho$. For example,  if $R$ is a
coset of a subgroup $H$ of finite index in $F$ then it follows
from Woess \cite{woess} that $$ \mu_0(R) = \frac{1}{|G:H|}, $$
while, obviously,  $\rho(H) = 1$ for the group $H$ of index $2$
consisting of all elements of even length.

 On the other hand, if $\lim_{k\rightarrow \infty} f_k(R)$
exists (hence is equal to $\rho(R)$) then $\mu_0(R)$ also exists
and $\mu_0(R) = \rho(R)$. In particular, if a set $R$ is
$\lambda$-measurable, then it is $\mu_0$-measurable, and
$\mu_0(R) = 0$.

\subsection{Asymptotic classification of subsets}
\label{se:1-5}

In this section we introduce  a classification of subsets $R$ in
$F$ according to the asymptotic behaviour of the functions
$\mu(R)$.

Let $\mu = \{\mu_s\}$ be the family of measures defined in
Section \ref{se:1-2}.
We start with a  global characterization of subsets of $F$.

Let $R$ be a subset of $F$. By its construction, the function
$\mu(R)$ is analytic on  $(0,1)$. The subset $R$ is called
{\it rational, algebraic, etc,  } with respect to $\mu$  if the
function $\mu(R)$ is  rational, algebraic, etc. We say that $R$ is
{\em smooth} if $\mu(R)$ can be analytically extended to a
neighborhood of $0$ and is regular at $0$.

\paragraph{Algebraic sets and context free languages.}
If the set $R$ is an (unambiguous) context free language then, by
a classical theorem of Chomsky and Schutzenberger \cite{Ch-Sch},
the generating function $\mu^*(R) = \sum n_kt^k$, and hence the
function $\mu(R)$, are algebraic functions of $s$. Moreover, if
$R$ is regular  then $\mu(R)$ is a rational function with
rational coefficients \cite{flajolet,stanley}.

 An important class of example of algebraic subsets is
provided by a theorem of Muller and Schupp \cite{muller-schupp}: A
normal subgroup $R \triangleleft F$ is a context free language if and only
if the factor group $F/R$ is free-by-finite. Notice that, for the
derived subgroup $R=[F,F]$ of the free group of rank $2$, the
measure $\mu(R)$ is not an algebraic function. Richard Sharp
kindly informed us that this follows from a remark on p.~127 of
his paper \cite{sharp}. See also Example~2.

It is well known that singular points of an algebraic function
are either poles
or branching points. Since $\mu(R)$ is bounded for
$s \in (0,1)$, this means that,
for a context-free set $R$, the function $\mu(R)$
has no singularity at $0$ or has a branching
point at $0$. After uniformisation, we can expand $\mu(R)$
 as a fractional power
series:
 $$ \mu(R) = m_0 + m_1 s^{1/n} + m_2s^{2/n} + \cdots . $$
If $R$ is regular, than we actually have the usual power series expansion:
 $$ \mu(R) = m_0 + m_1 s + m_2s^{2} + \cdots; $$
 in particular,  $\mu(R)$ can be analytically extended in the
 vicinity of $0$ and $R$ is smooth.

\paragraph{Linear approximation.}
 If the set $R$ is smooth then  the linear term in the expansion  of
$\mu(R)$ gives a linear approximation of  $\mu(R)$:
$$
\mu(R) = m_0 + m_1s + O(s^2).
$$
 Notice that, in this case,  $m_0 = \mu_0(R)$ is the Cesaro density of $R$.
 It can be shown (see Corollary~\ref{cor:Cesaro2} in Section~\ref{se:6})  that
 if $\mu_0(R) = 0$ then
$$ m_1 = \sum_{k=1}^\infty f_k(R) = \lambda(R). $$

On the other hand, even without assumption that $R$ is smooth,
 if $R$ is $\lambda$-measurable (that is, the series $\sum
f_k(R)$ converges), then, by Corollary~\ref{cor:Cesaro2},
 $$\mu_0(R) = 0 \ \  \hbox{ and }\ \ \lim_{s\rightarrow 0^+} \frac{\mu(s)}{s} =  \lambda(R).$$
This give us a good excuse to use for the limit
$$\mu_1 = \lim_{s\rightarrow 0^+} \frac{\mu(s)}{s},$$
if it exists, the same term {\em
frequency measure} as for $\lambda$. The function $\mu_1$ is an
  additive   measure on $F$ (though it is not $\sigma$-additive).

\paragraph{Asymptotic classification of sets.}
 Now we introduce a subtler classification of sets in $F$
 (which is based on the linear approximation of $\mu(R)$:
\begin{itemize}
\item {\em Thick} subsets: $\mu_0(R)$ exists, $\mu_0(R) > 0$ and
$$\mu(R) = \mu_0(R) + \alpha_0(s), \ \ where \ \ \lim_{s \rightarrow
0^+}\alpha_0(s) = 0.$$

\item {\em Sparse} subsets: $\mu_0(R) = 0$, $\mu_1(R)$ exists  and
$$\mu(R) = \mu_1(R)s + \alpha_1(s)\ \  where \ \ \lim_{s \rightarrow
0^+}\frac{\alpha_1(s)}{s}  = 0.$$

\item {\em Intermediate density} subsets: $\mu_0(R) = 0$
 but $\mu_1(R)$ does not exist.

\item {\em Singular\/} sets: $\mu_0(R)$ does not exist.

\end{itemize}

We put on record the following simple observation which
follows from discussions in Section~\ref{se:1-2}.

\begin{lemma}
Every $\lambda$-measurable set is sparse.
In particular, if $\gamma(R) < 1$ then $R$ is sparse.
\label{lm:lambda=sparse}
\end{lemma}

 We shall see in  Section~\ref{sec:regular} that, for the important
class of regular sets, the generating function is a rational
function and hence every regular set is either thick or sparse.

\subsection{Degrees of polynomial growth}
\label{se:1-6}

In this section we  introduce degrees of polynomial growth ``on
average'' for functions on the free group with respect to the
family of distributions $\mu$. In particular, it would produce
 hierarchies of the {\it average case complexity} of
various algorithms for infinite  groups,
which would make meaningful statements
like ``the algorithm works in cubic time   on average''.
A different approach to degrees of growth ``on average''
was suggested in \cite{measures}.

Let  $\mu = \{\mu_s\}$ be  the family of measures constructed in Section
\ref{se:1-2} and $\lambda$ be the frequency measure on $F$.  Let  $f: F_n
\longrightarrow {\mathbb{R}}$ be a non-negative real valued function.

 The average value $E_k = E_k(f)$ of the function $f$ on the sphere $S_k$ (with respect to
$\lambda$) is equal to:
 $$ E_k = \sum_{w \in S_k} f(w)\lambda(w) = \sum_{|w| = k}
 \frac{f(w)}{|S_k|}.$$

For every fixed stopping probability $s\in (0,1)$ we evaluate the
mean value $M_f(s)$ of the function $f$ with respect to $\mu_s$ as
\begin{eqnarray*}
M_f(s) &=& \sum_{w\in F} f(w)\mu_s(w)\\  & = & s\sum_{k=0}^\infty
E_k(1-s)^k.
\end{eqnarray*}
If for every $s \in (0,1)$ the value $M_f(s)$ is finite then the
function $M_f(s)$ is called the {\it mean value } of $f$ with
respect to the family of distributions $\mu$. The growth of the
function $M_f$ at $s = 0$ corresponds to the growth of the mean
values of $f$ with respect to the family  $\mu$ when the mean
length $L = \frac{1}{s} -1$ tends to infinity. Therefore, if we
rewrite $M_f(s)$ in the variable $L$:
  $$ M^*_f(L) = M_f\left(\frac{1}{L+1}\right),  \ \ L \in (0,\infty),$$
then the growth of $M^*_f$ at $\infty$ reflects the growth of the
initial function $f$ when the length of words tends to $\infty$.
This allows one to introduce the notion of the polynomial growth
of $f$ on average.

 Let $\nu: (1,\infty) \rightarrow {\mathbb{R}}$ be an arbitrary
 continuous probability
 density  on the interval $(1,\infty)$  and $\nu(x)dx$ the corresponding
 probabilistic measure.  We say that a non-negative  real
 valued function $f: F \rightarrow {\mathbb{R}}$ has a
 {\em polynomial growth of degree
  $d$ on average with respect to $\mu$}  and $\nu$
  if the function $M^*_f$  has
polynomial growth of degree   $d$ on average with respect to
  $\nu$, i.e., the
following improper integral converges at $\infty$:
 $$ \int_{0}^{\infty} \frac{M^*_f(x)}{x^{d}} \nu(x) dx.$$
  and $d\in {\mathbb N}$ is the minimal with this property.

 If $\eta(s)ds$ is the measure on $(0,1)$ obtained from $\nu(x)dx$
 by the change
 of variables $s = 1/(x +1)$, then this is the same as to say that
 $$ \int_{0}^{1} s^d M_f(s)\eta(s) ds$$
 converges at $0$.
 In most cases we can use the standard measure $ds$ on $(0,1)$.

Elementary results from analysis give the following  simple and useful test for
polynomial growth of functions on average.
\begin{lemma}
\label{Le:poly-test}
 Let $f: F \rightarrow {\mathbb{R}}$ be a non-negative real valued function on $F$.
 If the mean value function $M_f(s)$ is defined for  $s \in (0,1)$ and,  in
the vicinity of $0$, $$ M_f(s) = O(s^{-d}) $$ for some positive
integer $d$ then $f$ has polynomial growth of degree at most $d$
on average
 for any continuous probabilistic measure $\eta(s)ds$ on $(0,1)$.
\end{lemma}

Our definition of polynomial growth on average is justified by the
following simple observation.

\begin{lemma}
The function $f(w) = |w|^n$ has growth of degree $n$ on average.
\label{lm:poli-growth}
\end{lemma}

Proof. In view of Lemma~\ref{Le:poly-test}, it will suffice to
prove that, in the vicinity of $0$,
$$
M_f(s) = O(s^{-n}).
$$
We shall work with a larger function
$$g(w)=(|w|+1)(|w|+2)\cdots (|w|+n).$$
Without loss of generality, we can assume $g(1) = 0$. It is easy
to see that its mean
$$
M_g(s) = s \sum_{k=1}^\infty (k+1)(k+2)\cdots (k+n)(1-s)^k.
$$
After changing the variable, $z = 1-s$, it is enough to prove that
the function
$$M(z)=(1-z)\cdot \sum_{k=1}^\infty (k+1)(k+2)\cdots (k+n)z^k$$
has  a pole of degree at most $d$ at $z=1$. But it is very easy
to see that
$$
M(z) =(1-z)\frac{d^n}{dz^n}\left(\frac{z^{n+1}}{1-z}\right)
$$
has a pole of degree $n$ at $z=1$. \hfill $\square$

\medskip

The following lemma shows that our definition of growth is natural
in the sense that polynomial growth of averages $E_k(f)$ of the
function $f$ over the spheres $S_k$ implies polynomial growth of
the function $f$ on average in the sense of our definition.

\begin{lemma}
If\/ $F_k(f) \leqslant Ck^d$ for some constant $C$, then $f(w)$
has polynomial growth of degree at most $d$.
\end{lemma}

Proof. Immediately follows from the previous lemma. \hfill
$\square$

\subsection{Negligible sets}

Let $R\subset F$ and $\chi$ be the characteristic function of $R$.
We say that a set $R$ is {\em polynomially negligible} if
 for
every positive integer $d$, the polynomial function $|w|^d\chi(w)$
restricted to $R$ has growth of degree at most $0$ on average.

The purpose of this concept is that, in computations of degrees of growth on
average, we can use a `cut and paste' technique and ignore any
polynomial function of any degree with support
restricted to $R$.

\begin{theorem}
Let $R\subset F$  and $f_k= |R\cap S_k|/|S_k|$ be the relative
frequency of elements of $R$ in the sphere $S_k$. Assume that the
function
$$
\mu^*(R) = \frac{2m}{2m-1}\sum_{k=0}^\infty f_k(1-s)^k
$$ can be continued analytically to a
neighborhood of\/ $0$ and is regular at $0$. Then the set $R$ is
polynomially negligible. \label{th:negligible}
\end{theorem}

\paragraph{Proof.} We can assume without loss of generality that
$R\subseteq F \smallsetminus\{1\}$
hence $f_0 =0$. Following the same line of  argument as in
Lemma~\ref{lm:poli-growth},  we replace $|w|^n$ with the larger
function
$$g(w)=(|w|+1)(|w|+2)\cdots(|w|+n)$$
and set $g(1)=0$.
 Let $z = 1-s$ and $r(z) = \sum f_k z^k$.  Observe that $r(z)$ is analytic
 and regular in a neighbourhood of
$z = 1$. It is
enough to prove that
$$
G(z) = (1-z)\cdot\sum_{k=0}^\infty (k+1)(k+2)\cdots (k+n)f_k z^k
$$
is regular at $z=1$. But this is obvious because $$G(z) =
(1-z)\frac{d^n}{dz^n}\left(z^{n+1} \cdot r(z)\right)$$ is regular
at $z=1$. \hfill $\square$

\begin{corollary}
If the
relative growth rate $\gamma(R) < 1$ then $R$ is negligible.
\label{cor:negligible}
\end{corollary}

\paragraph{Proof.} Since the radius of convergence $r$ of the series
$\sum f_k(R) (1-s)^k$ is computed as
$r = 1/\gamma(R)$,
we see that the function $\sum f_k(R) (1-s)^k$ is analytic and
regular in the vicinity of $0$. Hence $R$ is negligible by
Theorem~\ref{th:negligible}.   \hfill $\square$

\section{Normal subgroups and cogrowth}
\label{sec:cogrowth}

\subsection{Non-recurrent and non-amenable factor groups}

 Let $G$ be a finitely generated group with an atomic
probability measure  $\nu: G \rightarrow [0,1]$. The measure $\nu$
is called {\it symmetric} if $\nu(g) = \nu(g^{-1})$ for all $g \in
G$. The {\it support} of $\nu$ is defined as $${\rm supp}(\nu)= \{g \in
G \mid \nu(g) \neq 0\}.$$ With a given measure $\nu$ on $G$ one
can associate a random walk $W_{\nu}$ on $G$ such that the
transition probability from $g$ to $h$ is equal to $\nu(g^{-1}h)$.
  A finitely generated group $G=F/R$ is called {\em recurrent}, if
it admits a symmetric atomic probability measure $\nu:G
\rightarrow [0,1]$, whose support generates $G$, and such that the
corresponding random walk $W_\nu$ on $G$ is recurrent. Recall,
that a random walk  is recurrent if it returns to $1$ infinitely
many times with probability $1$, i.e., the series
$$\sum_{n= 0}^{\infty} p^{(n)}(1),$$
where $p^{(n)}(1)$ is the probability for the walker to return to
$1$ in $n$ steps,  is divergent.   By a result of Varopoulos
\cite{varopoulos} based on Gromov's polynomial growth theorem
\cite{gromov}, a group is recurrent if and only if it is finite or
a finite extension of $\mathbb{Z}$ or ${\mathbb{Z}}^2$.
Grigorchuk gave in
\cite{grigorchuk} another useful characterization
of recurrent groups:  the group $G = F/R$ is recurrent if and only
if the series $$\sum_{k=0}^\infty \frac{n_k(R)}{(2m-1)^{k}}$$
diverges. Observe, that the latter is equivalent to the condition
that $\sum_{k=0}^\infty f_k(R)$ is divergent, i.e., $\lambda(R) =
\infty$.

\begin{theorem}
Let  $R$ be a normal subgroup in a free group $F$. If the
factor group $F/R$ is not recurrent then $R$ is sparse.
\label{th:recurrence}
\end{theorem}

\paragraph{Proof.} Let  $F$ be a free group of rank $m$, $R$
be a normal subgroup of $F$ such that $F/R$ is not recurrent. Let
$n_k = |R \cap S_k|$, and $f_k = n_k/|S_k|$. Since   $F/R$ is not
 recurrent, it is infinite and, by a result of Woess
\cite{woess}, the asymptotic density
$\rho(R) = \lim_{k\rightarrow \infty}  f_k$ exists
and equal to $0$. By the  criterion above  the series
$\sum_{k=0}^\infty f_k $  converges. Therefore it is  Abel
summable, i.e.,  there exists a limit
$$\lim_{s\rightarrow 0^+} \sum_{k=0}^{\infty} f_k(1-s)^k =
\sum_{k = 0}^{\infty} f_k$$ and
\begin{eqnarray*}
\mu(R) & = &  s \cdot  \sum_{k=0}^{\infty}
f_k(1-s)^k \\
& = & \mu_1 s +o(s),
\end{eqnarray*}
where $$\mu_1 = \sum_{k=0}^{\infty}
f_k.
$$
\hfill $\square$

\medskip

A classical criterion of amenability, due to Cohen \cite{cohen}
and Grigorchuk \cite{grigorchuk}  claims that a finitely
generated group $F/R$ is {\it amenable} if and only if the {\em
cogrowth coefficient\/} $\lim\sup (n_k(R))^{1/k}=
2m-1$. This immediately gives the following result.

\begin{theorem}
Let $R$ be a normal subgroup in $F$. If the factor group $F/R$ is
not amenable then $R$ is sparse and polynomially negligible.
\label{th:amenability}
\end{theorem}

\paragraph{Proof.}  We have mentioned in Section \ref{se:1-2}
that if $\gamma(R) < 1$ then $R$ is $\lambda$-measurable, i.e.,
the series $\sum_{k = 0}^{\infty}f_k(R)$ converges. In the same time,
\begin{eqnarray*}
\frac{\lim\sup \left(n_k(R)\right)^{1/k}}{2m-1} & = &
\lim\sup \left(\frac{n_k(R)}{(2m-1)^k}\right)^{1/k}\\
& = & \lim\sup \left(\frac{n_k(R)}{2m(2m-1)^{k-1}}\right)^{1/k}\\
& = & \lim\sup \left(f_k\right)^{1/k}\\
& = & \gamma(R)
\end{eqnarray*}
Hence if $\lim\sup (n_k(R))^{1/k} < 2m-1$ then $\gamma(R) < 1$ and $R$ is
sparse by Lemma~\ref{lm:lambda=sparse}.
By Corollary~\ref{cor:negligible}, $R$ is polynomially negligible and
the theorem follows.
\hfill $\square$
\medskip

It is worth mentioning a corollary from the proof: since
 the convergence radius of the generating function
$$N(R) = \sum_{k=0}^\infty n_kt^k,
$$
is $\left(\lim\sup \left(n_k(R)\right)^{1/k}\right)^{-1}$,
we have:
\begin{corollary}
The convergence radius of $N(R)(t)$ is $((2m-1)\gamma(R))^{-1}$.
\label{cor:radius}
\end{corollary}

\subsection{Return generating function}

Let $M = M(X)$ be the free monoid generated by $X \cup
X^{-1}$ and $$\eta: M(X) \rightarrow F(X)$$ the canonical
epimorphism of monoids which is induced  by the identity map on $X
\cup X^{-1}$. Analogously to free groups,  denote by $S_k(M)$ the
set (sphere) of all words in $M$ of length $k$.  Then, given a
normal subgroup $R \triangleleft F$, we can consider two generating
functions,
$$
N(t) = \sum n_k t^k \; \hbox{ and } \; N^*(t) = \sum n^*_kt^k,
$$
where $n_k = n_k(R) = |R \cap S_k|$ and $n^*_k = n^*_k(R) =
|\eta^{-1}(R) \cap S_k(M)|$ (the number of words of length $k$ in
$M$ which are mapped  into $R$ by $\eta$). The function $N^*(t)$
is called the {\em return generating function}.

The following formula links
 the functions $N(t)$ and $N^*(t)$ for a
normal subgroup $R \triangleleft F$:
\begin{equation}
\frac{N(t)}{1-t^2} = \frac{N^*\left(\frac{t}{1+(2m-1)t^2}\right)}{1+(2m-1)t^2}.
\label{eq:godsil}
\end{equation}
(In  \cite{bartholdi} Bartholdi proved a more general result (see
Section~\ref{sec:non-normal}),  attributing
Equation~(\ref{eq:godsil}) to  Godsil \cite[p.~72]{godsil}.)

\begin{quote}
{\small {\bf Example 1.} We shall use Equation~(\ref{eq:godsil})
for the computation of
the measure of the {\em co-diagonal subgroup} $D$ of $F$, that is, the kernel
of the homomorphism $F \longrightarrow {\mathbb Z}$ defined by mapping
all generators $x_i \in X$ to the generator $1$ of ${\mathbb Z}$.

It is easy to see that there are $\bnm{2k}{k}m^{2k}$ words of
length $2k$ in $M$ which are mapped  by $\eta$ into $D$. Indeed,
these words are $2k$-tuples of elements $x_i$ (and there are
$m^{2k}$ of them) with the exponents $\pm 1$ assigned to them such
that the sum of exponents is $0$; there are $\bnm{2k}{k}$
assignments of exponents. Since
$$
\sum_{k=0}^\infty \bnm{2k}{k}t^{2k} = \frac{1}{\sqrt{1-4t^2}}
$$
(see sequence A000984 of \cite{online}),
$$
N^*(t) = \sum_{k=0}^\infty \bnm{2k}{k}m^{2k}t^{2k} =
\frac{1}{\sqrt{1-4m^2t^2}},
$$
and
\begin{eqnarray*}
N(t) &=&
\frac{1-{t^2}}{\left(1+(2m-1) {t^2}\right) {\sqrt{1-\frac{{4m^2} {t^2}}{{{(1+(2m-1) {t^2})}^2}}}}}\\
&=& \sqrt{\frac{1-t^2}{1-(2m-1)^2t^2}}
\end{eqnarray*}
A close look at the zeroes of the denominator
in the expression for $N(t)$ tells us
that the convergence radius of $N(t)$ is $1/(2m-1)$,
and Corollary~\ref{cor:radius} reinterprets
this statement as $\gamma(D)=1$.
A direct computation shows that
$$
\mu(D) = \frac{(m-1)}{m\sqrt{2-\frac{2}{m}}}\cdot \sqrt{s} + o(\sqrt{s}).
$$
An analysis along the lines of
the proof of Theorem~\ref{th:negligible} shows
that $D$ is not polynomially negligible.
}
\end{quote}

Notice that $$p_k = \frac{n^*_k}{(2m)^k}$$ is the probability for a
simple random walk\footnote{This means that we move from a vertex
to any of $2m$ adjacent vertices with equal probabilities $1/2m$.}
on $\Gamma$ to return to the initial vertex $H$ after $k$ steps.
A considerable body of literature on random walks on groups contains
various information about the {\em return probability generating function}
$P(t) = \sum_{k=0}^\infty p_kt^k$ (which is a special instance of the
{\em Green function} of the random walk).
Since $N^*(t) = P(2mt)$ and the frequency generating function
$$F(t) = \sum_{k=0}^\infty f_kt^k = 1 + \sum_{k=1}^\infty \frac{n_k}{2m(2m-1)^{k-1}}t^k
$$ is related to $N(t)$ as
$$
F(t) = 1+ \frac{2m-1}{2m}\left(N\left(\frac{t}{2m-1}\right) -1\right),
$$
 we easily convert
(\ref{eq:godsil}) into the following formula:
\begin{equation}
F(t) = \frac{1}{2m}\cdot \frac{(2m-1)^2 - t^2}{(2m-1)+t^2}
\cdot P\left(\frac{2mt}{(2m-1)+t^2} \right) + \frac{1}{2m}.
\label{eq:return-frequency}
\end{equation}

%\subsection{More examples}
\begin{quote}
{\small {\bf Example 2.}  For a normal subgroup $N \triangleleft F$,
Equation~\ref{eq:return-frequency} reduces the question of algebraicity of the function $\mu(N)$ to that one for  the generating function $N^*(t) = \sum n^*_kt^k$ for the number of non-reduced words in $N$. Assume that $F$ has rank $2$ and take for $N$ the derived subgroup
$N = [F,F]$ of $F$. Then non-reduced words of length $k$ from the
free monoid $M$ correspond to simple random walks of length
$k$ on the factor group $F/N \simeq {\mathbb{Z}} \times  {\mathbb{Z}}$
which start and end at $0$, the probability of that event being $n^*_k/(2m)^k$.
We found ourselves in the classical realm of random walks on lattices.
A paper by Montrol \cite[p.~201]{montrol} (see also \cite{maassarani})
contains a closed formula for the return probabilities generating function for a simple random walk
on  ${\mathbb{Z}} \times  {\mathbb{Z}}$:
$$
P(t) = \frac{2}{\pi z} Q_{-1/2}\left(\frac{2-t^2}{t^2}\right),
$$
where $Q_{-1/2}(z)$ is a Legendre function of the second kind.
As shown in \cite[Equation~22 on p.~201]{montrol},
$$
P(t) \sim - \frac{1}{\pi t} \log\left(\frac{1-t^2}{t^2} \right)
\quad \hbox{as} \quad t \rightarrow 1.
$$
Since $\mu_s(N) = s F(1-s)$, after an easy calculation with Equation \ref{eq:return-frequency},
we see that the function $\mu(N)$ is not algebraic and
$$
\mu_s(N) \sim \frac{s\log s}{\pi} \quad \hbox{as} \quad s \rightarrow 0.
$$
In particular, $\mu_s(N)/s$
has logarithmic divergence at $s=0$, $\mu_1$ does not exists and $N=[F,F]$
is a subgroup of intermediate density.
}
\end{quote}

\subsection{Non-normal subgroups and random walks on regular graphs}

\label{sec:non-normal}
 In this section, we transfer
Theorem~\ref{th:amenability} from normal to arbitrary subgroups of
$F$.

If $H$ is a (not necessarily normal) subgroup of $F = F(X)$, the
set $F/H$ of right cosets gives rise to the Schreier graph of
$H$,  denoted by $\Gamma$, if we connect the cosets $Hy$ and
$Hyx$, $x \in X$, by a directed edge marked $x$. Every closed path
in $\Gamma$ from $H$ to $H$ gives  a word in  the free monoid
$M(X^{\pm 1})$ which represents an element from $H$, if we read
the edge label when we go along the edge, and its inverse, if we
go against the direction of the edge. Reduced words correspond to
paths without backtracking of edges.

Notice that $\Gamma$ is is a $2m$-regular graph, that is, every
its vertex has valency $2m$.

Denote by $n_k$ the number of closed paths without backtracking of
edges which start and end at the vertex $H$. Notice that $n_k =
n_k(H)$ is exactly the number of reduced words of length $k$ in
$H$. Also, denote by $b_k = b_k(H)$ the number of all paths
of length $k$  from $H$ to $H$, and let
$$
N(t) = \sum_0^{\infty} n_kt^k \quad\hbox{ and }\quad B(t) =
\sum_0^{\infty} b_kt^k
$$
be the corresponding generating functions.

 Formula (\ref{eq:godsil}) is a special case of the following
result valid for all regular graphs \cite{bartholdi}:
\begin{equation}
\frac{N(t)}{1-t^2} =
\frac{B\left(\frac{t}{1+(2m-1)t^2}\right)}{1+(2m-1)t^2}.
\label{eq:bartholdi}
\end{equation}
Denote by $$p_k = \frac{b_k}{(2m)^k}$$  the probability for a
simple random walk
on $\Gamma$ to return to the initial vertex $H$ after $k$ steps.
The quantity
$$
\nu = \lim\sup \sqrt[k]{p_k}
$$
is called the {\em spectral radius} of $\Gamma$. Obviously, $\nu
\leqslant 1$.

\begin{theorem}
If the coset graph $\Gamma$ of a subgroup $H < F$ has spectral
radius $\nu < 1$ then $H$ is sparse and polynomially negligible.
\label{th:nonnormal}
\end{theorem}

\paragraph{Proof.} Let $r'$ and $r''$ be the convergence radii
of the formal power series  $N(t)$  and $B(t)$, then, by the
well-known result from calculus,
$$ r' = (\lim\sup \sqrt[k]{n_k})^{-1} \quad\hbox{ and }\quad
r'' = (\lim\sup \sqrt[k]{b_k})^{-1}.$$ If $\nu < 1$ then,
obviously, $r''
> 1/(2m)$. The formula (\ref{eq:bartholdi}) relates the convergence radii
of $N(t)$ and $B(t)$ (see also \cite{northshield} where this
relation was developed earlier). It is easy to see that $r' >
1/(2m-1)$, hence
 for the relative growth rate of $H$ we have:
$\gamma(H) < 2m-1$, and by Corollary
\ref{cor:negligible}  $H$ is sparse and polynomially negligible.
\hfill $\square$

\medskip

 A similar technique with the use of results from
\cite[Chapter~7]{kitchens} proves the following theorem.

\begin{theorem} {\bf (T.~Smirnova-Nagnibeda, private communication)}
If\linebreak $H <F$ is a subgroup of infinite index then its asymptotic
density $$ \rho(H) = 0.
$$
In particular, it follows that $\mu_0(H) = 0$.

\end{theorem}

\subsection{\relax Preimages of quasiconvex
subgroups\\ of hyperbolic factor groups are negligible}

Recall that a finitely generated group $G$ is {\em
word-hyperbolic} if for any  (some) finite generating set $S$
of $G$ there is $\delta \geqslant 0$ such that all geodesic
triangles in in the Cayley graph $C(G,S)$ of $G$ with respect to
$S$ are $\delta$-thin, that is, each side is contained in the
closed $\delta$-neighbourhood of the union of the other two sides.
A subgroup $H$ of a word-hyperbolic group $G$ is {\em quasiconvex}
if for any  (some) generating set $S$ of $G$ there is
$\epsilon \geqslant 0$ such that every geodesic in $C(G,S)$  with
both endpoints in $H$ is contained in the $\epsilon$-neighbourhood
of $H$.

 In \cite[Theorem~1.2]{kapovich} I.~Kapovich proved  that the
coset graph $G/H$ of a quasiconvex subgroup $H$ of a hyperbolic
group $G$ has spectral radius $<1$. Now, as an application of
Theorem~\ref{th:nonnormal} we have the following result.

\begin{theorem}
Let $R\triangleleft F$ be a normal subgroup of a free group $F$ such that
$F/R$ is a non-elementary word-hyperbolic group, and $R \leqslant
H < F$ a subgroup of infinite index in $F$ such that $H/R$ is a
quasiconvex subgroup of $F/R$. Then $H$ is sparse and polynomially
negligible in $F$.
\end{theorem}

\section{Measure of a regular set}
\label{sec:regular}

\subsection{Regular Languages and finite automata}

In this section we show how to compute the measure $\mu(R)$
of a regular subset of the free group $F = F(X)$ of rank $m$.
Most of the results here are just a proper interpretation of some
well-known facts about regular sets. We refer to \cite{epstein}
for detailed discussion of regular subsets of $F$.

Recall that a {\em finite automaton} $\cal A$ is a finite
labelled  oriented graph (possibly with multiple edges and
loops). We refer to its vertices as {\em states\/}. Some of the
states are called {\em initial\/} states, some {\em accept\/}
states.  We assume further that every edge of the graph is
labelled by one of the  symbols  $x^{\pm 1}, x \in X$. A {\em
path in} {\cal A}  is a set of edges  $e_0, \ldots, e_l$ such
that, for each $i=1,\ldots, l$, the endpoint of $e_{i-1}$ is equal
to the starting point  of  $e_i$. Reading the labels on edges
along the  path  in the natural order, we get the label of the
path. The {\em language accepted by an automaton} ${\cal A}$ is
the set ${\cal L } = {\cal L(A)}$ of labels on paths from an
initial state to an accept state. An automaton is said to be {\em
deterministic} if, for any state, there is at most one arrow with
the given label exiting from the state. A {\em regular} set  is a
language accepted by a finite deterministic automaton. For every
 finite deterministic automaton ${\cal A}$  one can construct
a finite deterministic automaton ${\cal A^*}$ such that ${\cal
L(A) = L(A^*)}$  and where the sets of initial and accept states
are disjoint. It would be convenient for us to work only with
non-empty words, that is, elements in $F \smallsetminus \{1\}$.

We assemble here some (mostly well known) facts about regular
sets.

\begin{theorem} Let\/ $A$ and\/ $B$ are regular subsets in $F$.
\begin{itemize}
\item The sets $A \cup B$, $A \cap B$ and\/ $A\smallsetminus B$ are regular.
\item The prefix closure $\bar A$ of a regular set\/ $A$ is regular. Here, the\/ {\em prefix
closure} $\bar A$ is the set of all initial segments of all words
in $A$.
\item If\/ $C$ is a regular set in the free monoid\/ $M$ freely generated by $X \cup X^{-1}$
then its image $\bar C$ under the natural reduction homomorphism
$M \longrightarrow F$ is regular.
\item The  product $$AB = \{\, ab \mid a\in A, \, b \in B \,\}$$ and the set of inverses
$$A^{-1} =\{\,a^{-1} \mid a\in A\,\}$$ are regular.
\item Every finite subset in $F$ is regular.
\item If $\phi: F \longrightarrow F$ is an endomorphism then the set $\phi(A)$ is regular.
\end{itemize}
\end{theorem}

Now we will show how to compute the measure $\mu({\cal L})$
of a regular language accepted by  a finite deterministic
automaton ${\cal A}$. Recall, that the measure $\mu =\{\mu_s\}$
gives rise to a multiplicative measure $\mu^* = \{\mu_s^*\}$
$$
\mu_s^*(w) = \left(\frac{2m}{2m-1} \cdot \frac{1}{s} \right)\cdot
\mu_s(w),
$$
such that
$$
\mu_s^*(w) = t^{|w|}, \ \ \ where \  t =  \frac{1-s}{2m-1}.
$$

 By numbering the states by numbers $1,\ldots, n$, we can
associate with the automaton $\cal A$ its {\em adjacency matrix}
$A$ by taking  an $n\times n$ matrix and writing the number of
arrows from from state $i$ to state $j$ in the position $(i,j)$.
It is easy to see that the number of different paths of length $l$
from state $i$ to state $j$ is $(A^l)_{ij}$ and the measure of the
set of labels on these paths is $t^l(A^l)_{ij}$. Let $I$ and $J$
be the sets of initial and accept states. If we denote $T = tA$
then it follows that
$$ \mu_s^*({\cal L}) =   \sum_{i \in I, \,j \in J}((T)_{ij}+(T^2)_{ij}+
\cdots ).
$$
In particular, the series on the right converges for every given
$s$.  Denote by $B = T+T^2+\cdots$ the matrix with entries
from the ring of formal power series $\mathbb{R}[[t]]$, then,
obviously,
$$
\mu^*({\cal L}) =   \sum_{i \in I, \,j \in J}B_{ij}
$$
and
$$ B = (I_n-T)^{-1} - I_n.$$

We come to the following formula:
$$ \mu^*({\cal L}) = \sum_{i \in I, \,j \in J}((I_n-T)^{-1}-I_n)_{ij}.
$$
If we replace $\cal A$ by the  automaton ${\cal A^*}$  which
accepts the same language and where an initial state is never an
accept state, we can simplify the formula and write
$$
\mu^*({\cal L}) = \sum_{i \in I, \,j \in J}((I_n-T)^{-1})_{ij}.
$$

We have as a corollary the following result.
\begin{theorem}
The measure $\mu^*(R)$ {\rm (}and hence the probability measure
$\mu(R)${\rm )} of a regular subset of $F$ is a rational
function in \/ $t$  {\rm (}and hence in $s${\rm )} with rational
coefficients.
\end{theorem}

We can now apply this theorem to the  Cesaro density (see
Section~\ref{se:1-4}) and asymptotic classification of regular
sets (see Section~\ref{se:1-5}).

\begin{corollary}
The Cesaro density of a regular set is a rational number.
\end{corollary}

\begin{corollary}
Every regular set is either thick or sparse.
\label{cor:thick-or-sparse}
\end{corollary}

\subsection{Thick regular sets}

 We describe below thick regular sets.

 A {\em cone} $C=C(w)$ with the vertex $w$ is a set of all
elements in $F$ containing the given word $w$ as initial segment.
Obviously, cones  are regular sets.

 Let  $B = \{ u \mid |u| \leqslant k-1\}$ be the ball of radius $k-1$. Then $F
\smallsetminus B$ is the union of $|S_k|$ cones each of which has
the same measure as  the given cone $C = C(w)$ with $|w| = k$.
Hence
$$
\mu(C(w)) = \frac{1-\mu(B)}{|S_{|w|}|} = \frac{1}{|S_{|w|}|} +
O(s).
$$
In particular, a cone is a thick regular set. The following
theorem shows that every thick regular set involves a cone.

\begin{theorem}
Let\/ $R$ be a regular subset of\/ $F$. Then $R$ is thick if and
only if its prefix closure $\bar{R}$ contains a cone.
\label{th:thick}
\end{theorem}

\paragraph{Proof.}   Notice that if a regular set $R$
 is accepted by a finite deterministic automaton $\cal A$, then its
prefix closure $\bar R$ is accepted by the automaton $\bar {\cal
A}$ obtained from $\cal A$ by extending the set of accept states
 by adding all states which belong to a directed path in $\cal A$
from an initial state of $\cal A$ to an accept state of $\cal A$.

Since cones are thick sets, one direction of our theorem
immediately follows from the following lemma.

\begin{lemma}
Let $R\subset F$ be a regular set. Then $R$ is thick if and only
if its prefix closure $\bar R$ is thick.
\end{lemma}

\paragraph{Proof.}
Of course, if $R$ is thick then $\bar R$ is thick. To prove the
reverse, we use the obvious observation that the union of finitely
many of sparse sets is sparse  (if $\mu_0$-measurable). Let
$\cal A$ be a finite deterministic automaton which accepts $R$ and
$v_1,\ldots, v_n$ the accept states of $\bar {\cal A}$. Denote by
$R_i$ the subset of $\bar R$ accepted by the state $v_i$. Then
$\bar R = R_1\cup\cdots\cup R_n$ and one of the regular sets $R_i$
is thick. If $w$ is a label on a directed path from $v_i$ to an
accept state, say $v_j$, of $\cal A$ then
$$ R_i \circ w = \{x\circ w \mid x \in R_i\}
$$
is obviously a thick set and belongs to $R$. \hfill $\square$
\medskip

\begin{figure}[h]
 \setlength{\unitlength}{0.8mm}
\begin{picture}(100,70)(-15,-10)

\thicklines \put(55,19){$\Longrightarrow$}

\put(0,0){\circle{6}} \put(0,20){\circle{6}}
\put(0,40){\circle{6}} \put(20,20){\circle{6}}
\put(40,10){\circle{6}} \put(40,30){\circle{6}}

\put(80,0){\circle{6}} \put(80,20){\circle{6}}
\put(80,40){\circle{6}} \put(100,10){\circle{6}}
\put(100,30){\circle{6}} \put(120,10){\circle{6}}
\put(120,30){\circle{6}}

\thinlines \put(3,3){\vector(1,1){14}}
\put(4,20){\vector(1,0){12}} \put(3,37){\vector(1,-1){14}}
\put(24,18){\vector(2,-1){12}} \put(24,22){\vector(2,1){12}}

\put(84,2){\vector(2,1){12}}
\put(84,22){\vector(2,1){12}}
\put(84,38){\vector(2,-1){12}}

\put(104,10){\vector(1,0){12}}
\put(103,13){\vector(1,1){14}}
\put(104,30){\vector(1,0){12}}
\put(103,27){\vector(1,-1){14}}

\put(19,19){\scriptsize {$A$}} \put(98,9){\scriptsize {$A''$}}
\put(98,29){\scriptsize {$A'$}}

\put(8,10){\small {$b$}} \put(8,21){\small {$a$}}
\put(10,31){\small {$a$}} \put(29,16){\small {$d$}}
\put(28,26){\small {$c$}}

\put(88,5.5){\small {$b$}} \put(88,26){\small {$a$}}
\put(89,36){\small {$a$}} \put(109,6){\small {$d$}}
\put(114,17.5){\small {$d$}} \put(111.5,24){\small {$c$}}
\put(109,31){\small {$c$}}

\end{picture}

\caption{{\small \sl Splitting the states of the automaton $\cal
A$.}} \label{RIS.Y}
\end{figure}
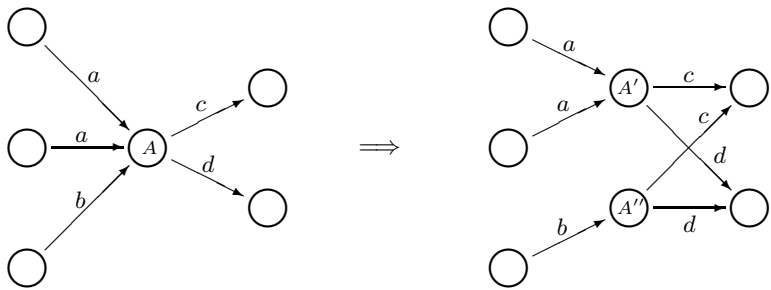

Now we can assume that the set $R$ is thick. Since the union of
finitely many sparse sets is sparse, we can assume without loss
of generality that a finite deterministic automaton $\cal A$ for
$R$ has only one initial state $I$ and one accept state $Z$. We
have to remember that our automaton accepts only reduced words.
Therefore $\cal A$ can be rewritten in the form where
\begin{itemize}
\item[(a)] For any state $A$ of $\cal A$,
all arrows which enter $A$ have the same label $a \in X \cup
X^{-1}$ and arrows exiting from $A$ cannot have label $a^{-1}$
(this can be achieved by splitting the states of $\cal A$ in the
way shown on Figure~\ref{RIS.Y}.) We shall say in this situation
that $A$ {\em has type} $a$.

\item[(b)] For every state $A$ of $\cal A$, there is a directed path
from $A$ to the accept state $Z$.

\item[(c)] In addition, it is easy to arrange that there are no arrows entering
the initial state $I$.
\end{itemize}

This means, in particular, that there are at most $2m$ arrows
exiting from the initial state $I$, and at most $2m-1$ arrows
exiting from any other state. We can assign frequencies $1/2m$ to
arrows exiting from $I$ and frequencies $1/(2m-1)$ to arrows
exiting from other states. Now, for a word $w$ accepted by $\cal
A$, its relative frequency
$$
\lambda(w) = \frac{1}{2m(2m-1)^{|w|-1}}
$$
is the product of frequencies of arrows in a directed path from
the initial state $I$ to the accept state $Z$ which correspond to
the word $w$. We aim at proving the following statement from
which our theorem immediately follows by virtue of
Lemma~\ref{lm:lambda=sparse}:

\begin{quote}
{\em  If $\bar R$ contains no cone then it is
$\lambda$-measurable, that is,
$$\lambda(R) = \sum_{w\in R} \lambda(w)$$ is finite. }
\end{quote}

For that purpose form the automaton ${\cal A}_1$ obtained from
$\cal A$ by removing all arrows exiting from $Z$; we take $I$ and
$Z$ for its initial and accept states, correspondingly. Consider
also the automaton ${\cal A}_2$ formed by all states accessible
from the state $Z$, with the same arrows between them as in $\cal
A$; we take $Z$ for the both initial and accept states.

  We assign to arrows in ${\cal A}_1$ and ${\cal A}_2$ the
same frequencies as to arrows in $\cal A$. Since $I$ does not
belong to ${\cal A}_2$, all arrows in ${\cal A}_2$ have
frequencies $1/(2m-1)$. If now $R_1$ and $R_2$ are languages
accepted by ${\cal A}_1$ and ${\cal A}_2$ then, obviously, $R =
R_1 \circ R_2$. Moreover, if $u \in R_1$ and $v \in R_2$ then the
word $uv$ is reduced and
$$
\lambda(uv) = \lambda(u)\lambda^*(v).
$$
Since the presentation of $R$ in the form $R = R_1 \circ R_2$ is
unambiguous, it follows that
$$
\lambda(R) = \lambda(R_1)\lambda^*(R_2).
$$
Transform the automaton ${\cal A}_2$ further by splitting the
state $Z$ into separate initial state $Z_1$ (with no arrows
entering it, and those arrows which exited from $Z$ now exiting
from $Z_1$), and the accept state $Z_2$ (with no arrows exiting
from it, and those arrows which entered $Z$ now entering $Z_1$).
If $R_3$ is the language accepted by the new automaton ${\cal
A}_3$, then, obviously,
$$
R_2 = R_3 \cup (R_3 \circ R_3) \cup (R_3 \circ R_3 \circ R_3) \cup
\cdots
$$
and
$$
\lambda^*(R_2) \leqslant \lambda^*(R_3) + \lambda^*(R_3)^2+
\lambda^*(R_3)^3 + \cdots .
$$

Assume that $\bar R$ contains no cone. Then the both subsets $\bar
R_1$ and $\bar R_2$ contains no cone.

Let us look first at ${\cal A}_2$. Assume that, for every state
$A$ of ${\cal A}_2$ of type $a \in X \cup X^{-1}$, every possible
label from $X \cup X^{-1} \smallsetminus \{a\}$ is present on  one
of the   arrows exiting from $A$. Then it is easy to see that
$\bar R_2$ contains a cone. Therefore we can assume that, for
some state $A$, there are less than $2m-1$ arrows exiting from
$A$. If we  now look at the automaton ${\cal A}_3$, it becomes
obvious that $\lambda^*(R_3) < 1$. To see this formally, we can
consider a Markov  chain $\cal M$  whose states are the states of
${\cal A}_3$ together with a additional dead state $D$ (see
\cite{ksk} for background material on Markov chains). We set the
transition probabilities from $Z_2$ to $Z_2$ and from $D$ to $D$
being equal $1$. Every arrow in ${\cal A}_3$ corresponds to a
transition in $\cal M$ with the transition probability
$1/(2m-1)$. If at some state $A$ of ${\cal A}_2$ there is no
arrow labelled $b \in X\cup X^{-1}$ exiting from $A$, we make in
$\cal M$ a transition from $A$ to $D$ with the transition
probability $1/(2m-1)$. The probability distribution on $\cal M$
concentrated at the initial state $Z_1$, converges to the steady
state $P$ which is zero everywhere with the exception of the two
dead states $Z_2$ and $D$. Since $P(D) \ne 1$, $P(Z_2) < 1$. But,
obviously,  $P(Z_2) = \lambda^*(R_3)$.

Now the summation of the geometric progression for
$\lambda^*(R_2)$ shows that $\lambda^*(R_2) < \infty$.

An analogous argument for ${\cal A}_1$ shows that $\lambda(R_1) <
\infty$. Therefore $\lambda(R) < \infty$.
  \hfill $\square$ $\square$

\subsection{Measures of finitely generated subgroups}

Let $F = F(X)$ be the free group with basis $X$. It is well known
that finitely generated subgroups in $F$ are regular sets;  the
most suitable for our purpose exposition of this and similar
results can be found in \cite{km}.

Let $\mu^*$ be the  adjusted multiplicative measure on $F$. Let
$H$  be a subgroup of $F$ generated by elements $h_1,\ldots,
h_k$. We shall slightly modify the  arguments  from the
previous section to produce a somewhat more practical procedure
for computing the measure $\mu^*(H)$. In particular, $\mu^*(H)$
will be expressed as a rational function of measures $\mu^*(w_i)$
of certain words $w_1, \ldots, w_s$ which do not depend on choice
of generators $h_1, \ldots, h_m$ in $H$ although can be easily
computed from them.

Let $\Gamma$ be the core subgroup graph of $H$ in sense of
\cite{km}. Notice that $\Gamma$ does not depend on a
particular choice of generators of $H$.  We mark on $\Gamma$ the
initial vertex $1$ and those vertices which have degree at least
$3$. This new vertex set $V^*$ can be turned into a digraph
$\Gamma^*$ with edges labelled by freely reduced words from
$F$. To do so, we define edges of $\Gamma^*$ to be  reduced paths
in $\Gamma$ which start and end at vertices in $V^*$ and  do not
pass through any other vertex from $V^*$. The label of the path
becomes the label of the corresponding edge in $\Gamma^*$. We call
$\Gamma^*$ the {\em consolidated subgroup graph} of $H$.

Now it is easy to see that, since $\Gamma$ is folded, a reduced
path in $\Gamma^*$, viewed as a path in $\Gamma$, is also reduced.
Every element  $h\in H$ is the   label of a reduced path in
$\Gamma$ from $1$ to $1$, as well as, the label $w_1\circ \ldots
\circ w_l$ of the corresponding reduced path in $\Gamma^*$ from
$1$ to $1$. It follows that $\mu^*(h) = \mu^*(w_1)\cdots
\mu^*(w_l)$.

 Our description of the  matrix method of computing the
measure of a finitely generated subgroup will be illustrated by
the following example, which we do in parallel with the formal
discussion.

\begin{quote}
{\small {\bf Example 3.} Let $C$ be  a subgroup  generated by a
single element $c$. Obviously $c$ can be presented in the form $c
= uvu^{-1}$ without cancellations between the words $u$, $v$ and
$u^{-1}$. The consolidated subgroup graph $\Gamma^*(C)$ of $C$ has
the form }
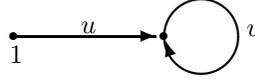
\begin{center}
\setlength{\unitlength}{.5mm}
\begin{picture}(100,40)(-20,0)
\thicklines
\put(19,18){$1$}
\put(20,25){\circle*{2}}
\put(20,25){\vector(1,0){38}}
\put(38,26){$u$}
\put(60,25){\circle*{2}}
%\put(63,23){$2$}
\put(60.5,24){\vector(-1,3){0}}
\put(70,25){\circle{20}}
\put(82,25){$v$}
\end{picture}
\end{center}
\end{quote}

We start with the consolidated subgroup graph $\Gamma^*$ of $H$.
If $e$ is an edge in $\Gamma^*$, we denote its label by
$\lambda(e)$. As usually, we use the convention that for every
edge $e$ we also have an edge with the opposite direction and the
inverse label $\lambda(e)^{-1}$,
\begin{quote}
{\small so, in our example, $\Gamma^*(C)$ is a digraph with $4$
directed edges. }
\end{quote}

The process of writing non-trivial random words from $H$ can be
described by the automaton $T$  which consists
of one state for each directed edge of the digraph $\Gamma^*$
plus one initial state. Every  directed edge $e$ of $\Gamma^*$ is
interpreted as the state ``we wrote the word $\lambda(e)$ of
$e$'', the initial state is ``we wrote the empty word''.  If the
origin of edge $f$ is the terminus of edge $e$, we say that there
is a transition from the state ``we wrote the word $\lambda(e)$
to the state ``we wrote the word $\lambda(f)$'', and we assign to
this transition measure $\mu^*(\lambda(f))$. To a directed edge
$e$ which exits from the initial vertex $1$ of $\Gamma^*$, we
assign the transition from the state ``we wrote the empty word''
to the state ``we wrote the word $\lambda(e)$'' with measure
$\mu^*(\lambda(e))$. Finally, the accept states of our automaton
 $T$ correspond to directed edges of $\Gamma^*$ whose
terminuses are the initial vertex $1$.

 We label
 the states of $T$ by consecutive numbers $1,2,3,\ldots$, so
that the initial state ``we wrote the empty word'' has label $1$.
Now the transition measures of automaton $T$
form a matrix which we denote $A$.

\begin{quote}
{\small In our Example 2, the states are:
\begin{itemize}
\item[1] We wrote an empty word (the initial state);
\item[2] we wrote $u$;
\item[3] we wrote $v$;
\item[4] we wrote $v^{-1}$;
\item[5] we wrote $u^{-1}$ (the accept state).
\end{itemize}
and the transition matrix is
$$
 A =
\left(\begin{array}{ccccc}
    0&p&0&0&0 \\
    0&0&q&q&0 \\
    0&0&q&0&p \\
    0&0&0&q&p \\
    0&0&0&0&0
\end{array}\right),
$$
where $p = \mu^*(u)$ and $q = \mu^*(v)$. }
\end{quote}

Let $j_1,\ldots, j_k$ be the accept states of $T$.
The multiplicativity of the adjusted measure $\mu^*$ allows to use
the same matrix technique as in computations on Markov chains,
and the measure of the set of words from $H$ which can be
obtained by $l$ moves becomes the sum
$$
\left(A^l\right)_{1j_1}+\cdots + \left(A^l\right)_{1j_k} $$ of
the matrix elements of the matrix $A^l$ which correspond to
transition from the initial to an accept state. Notice that we
produce only non-trivial elements of $H$. Hence the set
$H\smallsetminus \{1\}$ of non-trivial elements in $H$ has the
measure
$$
\mu^*(H\smallsetminus \{1\}) = \sum_{i=1}^k
\left(A+A^2+\cdots+A^n+\cdots\right)_{1j_i}.
$$
Denote $B = A+A^2+\cdots$, then, obviously,
$$
B = (E-A)^{-1} - E,
$$
where $E$ is the identity matrix. Since, by our construction, the
initial state is never an accept state, the matrix elements
$(B)_{ij_1},\ldots,(B)_{1j_k}$ do not lie on the diagonal and
therefore $(B)_{1j_i} = ((E-A)^{-1})_{1j_i}$ for all
$i=1,\ldots,k$. Hence
\begin{equation}
\mu^*(H\smallsetminus \{1\}) = \sum_{i=1}^k
\left((E-A)^{-1}\right)_{1j_i}.
\label{eq:main}
\end{equation}

Since the elements of the inverse matrix $(E-A)^{-1}$ are
rational functions of matrix elements of the matrix $A$,
we proved the following theorem.

\begin{theorem}
If $\mu$ is the multiplicative measure on $F$, the measure
$\mu(H)$ of a finitely generated subgroup $H$ of $F$ is a
rational function of the measures of labels on the consolidated
subgroup graph of $H$.

In particular, $\mu(H)$ is a rational function of $s$.
\end{theorem}

\begin{quote}
{\small In Example 3, a direct computation with
\textsc{Mathematica} shows that
$$
(E-A)^{-1} = \left(\begin{array}{ccccc}
    1&p&\frac{p q}{1-q}&\frac{p q}{1-q}&\frac{2 {p^2} q}{1-q} \\
    0&1&\frac{q}{1-q}&\frac{q}{1-q}&\frac{2 p q}{1-q} \\
    0&0&\frac{1}{1-q}&0&\frac{p}{1-q} \\
    0&0&0&\frac{1}{1-q}&\frac{p}{1-q} \\
    0&0&0&0&1
\end{array}\right)
$$
and
\begin{eqnarray*}
\mu^*(C\smallsetminus \{1\}) &=&
\frac{2p^2q}{1-q}\\
&=& \frac{2t^{2|u|}\cdot t^{|v|}}{1-t^{|v|}}.
\end{eqnarray*}
Hence
\begin{eqnarray*}
\mu(C\smallsetminus \{1\})
&=& \frac{2m-1}{2m}s\cdot \frac{2(\frac{1-s}{2m-1})^{2|u|+|v|}}{1-(\frac{1-s}{2m-1})^{|v|}}\\
&=& \frac{2}{2m(2m-1)^{2|u|+|v|-1} \cdot (1-(\frac{1}{2m-1})^{|v|})}\cdot s  +o(s).
\end{eqnarray*}}
\end{quote}

\subsection{Normal subgroups of finite index}

Kouksov \cite{kouksov} proved that a normal subgroup $N \triangleleft F$ has a rational
cogrowth function $f(t)=\sum n_kt^k$ if and only if the index $|F:N|$ is finite.
In that case $N$ is finitely generated, and its subgroup graph $\Gamma$ is the
Cayley graph of $G$. In notation of (\ref{eq:godsil}), the generating function $N^*(t)$ for the number of non-reduced words in $N$ has a very beautiful form found by
Quenell \cite{quenell}:
$$
N^*(t) = \frac{1}{|G:N|} \sum \frac{1}{1-\lambda_i t},
$$
where $\lambda_i$ are the eigenvalues of the adjacency matrix of $\Gamma$.

\section{Context free languages}

Combinatorial analysis of context free languages is a well
established area of combinatorics with powerful tools for
manipulating generating functions of languages; a very good
exposition of the theory, pioneered by Chomski and Schutzenberger
\cite{Ch-Sch}, can be found in \cite{flajolet,stanley}. Here we
give only a small example of use of this machinery, motivated by
applications of our methods to study of complexity of algorithms
on amalgamated products of free groups \cite{amalgam}. We do not
give rather technical and lengthy definitions related to
context-free languages which can be found in \cite{rayward-smith}
(see also \cite{woess2} for a compact formal definition).

Let $M$ be the free monoid generated by $X \cup X^{-1}$.
We call two subsets $R, S \subset M$ {\em isobaric} if, for every $k$, they contain equal number of words of
length $k$,
that is, if the have the same generating function.

\begin{quote}
{\small {\bf Example 4.} Let $X = Y \sqcup Z$, where $|Y|= l
\geqslant 1$ and $|Z|= n \geqslant 1$. We shall find the measure
of the set
$$
R = \bigcup_{g\in F}{F(Y)^\sharp}^g
$$
of all elements in $F$ conjugate to non-identity elements in $F(Y)$.
Here, as usually, we
denote
by $F(Y)^\sharp$ the set of non-identity elements of $F(Y)$.
Obviously, we can decompose
$$
R = {F(Y)^\sharp} \cup \bigcup F(Y)^h,
$$
where the union is taken over all elements in
$F$ which start with letters in $Z\cup Z^{-1}$.
The generating function for $F(Y)^\sharp$ is obvious:
\begin{eqnarray*}
f(t) & = &  \sum_{k=1}^\infty 2l(2l-1)^{k-1} t^k\\
& = & 2l\cdot t(1- (2l-1)t)^{-1}.
\end{eqnarray*}

Denote $L_1 = F(Y)^\sharp$ and let
$L_2$ be the cone of words in $F$ which start from symbols in $Z \cup Z^{-1}$.
It is easy to see that, in the free monoid $M(X)$,
the language $L_2^{-1} \circ L_1 \circ L_2 $ is isobaric to
to the language $L_1 \circ L_2^2$, where
$L_2^2  = \{f\circ f \mid f \in L_2 \}$ and $\circ$ denotes
formal product in $M(X)$ without cancellation.
The generating function of $L_2$ is
$$
g(t) = \sum_{k=1}^\infty 2n(2m-1)^{k-1}t^k = 2nt (1-(2m-1)t)^{-1},
$$
and the generating function for $L_2^2$ is $g(t^2)$.
According to the standard rules of computation of generating functions for
context-free languages \cite{stanley},
 the generating function for $L_1 \cup L_1\circ L_2^2$ is
$$
f(t) + f(t)g(t^2) = 2l\cdot t(1- (2l-1)t)^{-1}(1+2nt^2(1-(2m-1)t^2)^{-1})
$$
and therefore
{\scriptsize
\begin{eqnarray*}
\mu(R) & = & \frac{2m-1}{2m}s \cdot 2l\cdot \frac{1-s}{2m-1}\left(1- (2l-1)\cdot \frac{1-s}{2m-1}\right)^{-1}\\
&&\qquad \times \left(1+2n\left(\frac{1-s}{2m-1}\right)^2\left(1-(2m-1)\left(\frac{1-s}{2m-1}\right)^2\right)^{-1}\right)\\
&=& \frac{l}{m}\left(1- \frac{2l-1}{2m-1}\right)^{-1}\left(1+\frac{2n}{(2m-1)^2}\right)\left(1-\frac{1}{2m-1}\right)^{-1}\cdot s + o(s). \end{eqnarray*}
}
}
\end{quote}

Notice that by \cite{muller-schupp}, a normal subgroup $N \triangleleft F$
is context-free if and only if the factor group $F/N$ is
free-by-finite.

\section{Addendum: A Tauberian Theorem by Hardy and Littlewood}
\label{se:6}

We found ourselves in the context where generalised summation
methods for series are essential.

\begin{quote}
{\small {\bf Example 5.} Consider the
subgroup $H$ of index $2$ in $F$ which consists of all words of
even length in $F$. Let $n_k$ be the number of elements of length
$k$ in $H$ and
$$f_k = \frac{n_k}{|S_k|}
$$
be the relative frequency of elements of length $k$ from $H$
among all elements of length $k$ in $F$. Obviously,
$$
f_k = \left\{\begin{array}{rl} 1 & \hbox{ if } k \hbox{ is even}\\
0 & \hbox{ if } k \hbox{ is odd}
\end{array} \right. .
$$
One can easily see that
\begin{eqnarray*}
\mu(H) & = & s+ s \sum_{k=1}^\infty
(1-s)^{2k}\\
& = & s+\frac{s(1-s)^2}{1-(1-s)^2}\\
%& = & s+\frac{(1-s)^2}{2-s}\\
& = & \frac{1}{2} +\frac{s}{4} + \frac{s^2}{8} + \frac{s^3}{16} +
\cdots
\end{eqnarray*}
When $s \longrightarrow 0^+$, $\mu(H) \longrightarrow 1/2$. }
\end{quote}

To explain the rather expected appearance of $1/2$ as the `limit
probability' of the subgroup $H$ in Example 5, we need to invoke
one of the so-called Tauberian theorems by Hardy and Littlewood.

\begin{theorem} {\rm \cite[Theorems  94]{hardy}}
Let\/ $\{a_n\}$ be a  sequence of real numbers such that the sequence of partial sums
$$
S_n = a_0+\cdots + a_n
$$
is bounded from below.
Assume also that the limit\/ {\rm (}the Abelian sum of $\{a_n\}${\rm )}
$$
\lim_{x \rightarrow 1^-} \sum_{k=0}^\infty a_k x^k
$$
exists and equals $S$. Then the sequence $\{a_n\}$ is Cesaro
summable in the sense that the limit
$$
\lim_{n\rightarrow \infty} \frac{S_0+S_1+\cdots +S_n}{n+1}
$$
exists and equal $S$.
\label{th:Hardy}
\end{theorem}

\begin{corollary}
Assume that the sequence $\{f_n\}$ of non-negative real numbers
is bounded and the sum $\sum_{k=0}^\infty f_k (1-s)^k$ converges
for all\/ $0< s < 1$. Assume, in addition, that
 there exists the limit
$$
\lim_{s\rightarrow 0^+} s \sum_{k=0}^\infty f_k (1-s)^k = \mu_0.
$$
Then
$$\mu_0 = \lim_{k\rightarrow \infty} \frac{f_0
+f_1+\cdots+f_k}{k+1}
$$
\label{cor:Cesaro1}
\end{corollary}

\paragraph{Proof.} Set $t = 1-s$ and  rewrite
\begin{eqnarray*}
s \sum_{k=0}^\infty f_k (1-s)^k & = & (1-t)\sum_{k=0}^\infty f_k
t^k \\
& = & f_0 + (f_1-f_0)t + (f_2-f_1)t^2+\cdots ,
\end{eqnarray*}
and the series on the right converges for all $0 < t < 1$.
Moreover,
$$
\lim_{t \rightarrow 1^-} \left(f_0 + \sum_{k=1}^\infty (f_k
-f_{k-1})t^k\right) = \mu_0.
$$
Since the partial sums
$$f_0+(f_1-f_0)+\cdots+
(f_k-f_{k-1})= f_k$$
are bounded from below,
 the previous theorem yields
\begin{eqnarray*}
\mu_0 & = & \lim_{k\rightarrow \infty}
\frac{[f_0]+[f_0+(f_1-f_0)]+\cdots +[f_0+(f_1-f_0)+\cdots+
(f_k-f_{k-1})]}{k+1}\\
& = & \lim_{k\rightarrow \infty} \frac{f_0+f_1+\cdots+ f_k}{k+1} .
\end{eqnarray*}
 \hfill $\square$

 \medskip

Corollary~\ref{cor:Cesaro1} explains, in particular, that for our
subgroup $H$ of index $2$ in $F$, $\mu(H) = \mu_0 + o(1)$,
where
$$ \mu_0 = \lim_{n\rightarrow \infty}
\frac{1+0+\cdots +1+0}{2n} = \frac{1}{2}.
$$

 \medskip

The following corollary is an easy consequence of
Corollary~\ref{cor:Cesaro1} and Theorem~\ref{th:Hardy}.

\begin{corollary}
Assume that the sequence $\{f_k\}$ of non-negative real numbers
is bounded, the sum $\sum_{k=0}^\infty f_k (1-s)^k$ is convergent
for all\/ $0 < s < 1$ and
 the function
 $$
 \mu(s) = s\sum_{k=0}^\infty f_k (1-s)^k
 $$
 has the limit
 $$
 \lim_{s\rightarrow 0^+} \mu(s) = \mu_0.
 $$
Then
\begin{itemize}
\item[{\rm (a)}] $\mu_0$ is the Cesaro limit
$$
 \mu_0 = \lim_{k\rightarrow \infty} \frac{f_0
+f_1+\cdots+f_k}{k+1}.
$$
\item[{\rm (b)}]If $\mu_0 = 0$ and the limit
$$ \mu_1=\lim_{s\rightarrow 0^+}\frac{\mu(s)}{s}
$$
exists then the sum\/ $\sum f_k$ is convergent and
$$
\mu_1 = \sum_{k=1}^\infty f_k.
$$
\item[{\rm (c)}] If the series $\sum f_k$ converges,
 then $\mu_0$ and $\mu_1$ exist, $\mu_0 =
0$ and
$$\mu_1 = \sum_{k=0}^\infty f_k.$$

\item[{\rm (d)}] In particular, if the function $\mu(s)$ is analytic in
the vicinity of\/ $0$ and regular at $s = 0$, then, in the power series expansion at $s = 0$,
$$
 \mu_0 = \lim_{k\rightarrow \infty} \frac{f_0
+f_1+\cdots+f_k}{k+1},
$$
and, if $\mu_0 = 0$, the next coefficient is given by
$$
\mu_1 = \sum_{k=1}^\infty f_k.
$$
\end{itemize}  \label{cor:Cesaro2}
\end{corollary}

\paragraph{Proof.} (a) directly follows from \ref{cor:Cesaro1}.

For a proof of (b), notice that by Theorem~\ref{th:Hardy}
$$
\mu_1 = \raisebox{.5ex}{$^C$}\!\sum_{k=1}^\infty f_k,
$$
where the sum $\raisebox{.5ex}{$^C$}\!\sum f_k$ is understood in
the sense of the Cesaro limit of the partial sums $S_k=
f_1+\cdots+f_k$:
$$
\raisebox{.5ex}{$^C$}\!\sum_{k=1}^\infty f_k = \lim_{n\rightarrow
\infty}\frac{1}{n}(S_1+\cdots +S_n).
$$
But for non-negative series, Cesaro summability is equivalent
to the ordinary convergence, which yields the result.

For (c), assume that the series $\sum f_k$ converges.
By Abel's theorem on continuity of sums of power series, the
function $f(s)$ is continuous on the interval $[0,1)$ and hence
$$
\mu_1 = \lim_{s\rightarrow 0^+} f(s) = \sum_{k=0}^\infty f_k
$$
and
$$
\mu_0 = \lim_{s\rightarrow 0^+} sf(s) = 0.
$$
(d) is an immediate corollary of (a)- and (b).
\hfill $\square$

\subsection*{Acknowledgements}

The authors thank Ilya Kapovich, Ziad Maassarani, Guennadi Noskov, Mark Policott, Richard Sharp
  and Tatiana Smirnova-Nagnibeda
for useful discussions.

\bigskip

\normalsize

\vfill

 \noindent \textsf{Alexandre V. Borovik, Department of
Mathematics, UMIST, PO Box 88,\linebreak Manchester M60 1QD,
United Kingdom}

\noindent {\tt borovik@umist.ac.uk}

\noindent {\tt http://www.ma.umist.ac.uk/avb/}

\medskip

\noindent \textsf{Alexei G. Myasnikov,
 Department of Mathematics, The City  College  of New York, New York,
NY 10031, USA}

\noindent  {\tt alexeim@att.net}

\noindent {\tt http://home.att.net/\~\,alexeim/index.htm}

\medskip
\noindent \textsf{Vladimir N. Remeslennikov, Omsk Branch of
Mathematical Institute SB RAS,\linebreak 13 Pevtsova Street, Omsk
644099, Russia}

\noindent {\tt remesl@iitam.omsk.net.ru}

\end{document}